\documentclass{amsart}
\usepackage{amssymb}

\newtheorem{thm}{Theorem}[section]

\newtheorem{lem}[thm]{Lemma}

{}
\theoremstyle{remark}
\newtheorem{rem}{Remark}

\newcommand{\R}{\mathbb{R}}

\def \eqskip { \vspace*{2mm} }
\def \disps {\displaystyle}

\newcommand{\ka}{{\kappa}}

\newcommand{\dist}{\operatorname{dist}}

\newcommand{\ds}{\displaystyle}

\newcommand{\formsk}{\vspace*{1mm}}

\newcommand{\intom}{\int_{\Omega}}
\newcommand{\intdo}{\int_{\partial\Omega}}

\newcommand{\vdo}{|\partial\Omega|}
\newcommand{\vo}{|\Omega|}
\newcommand{\va}{\mathcal{B}}
\newcommand{\sa}{\mathcal{S}}

\begin{document}

\title[Characterization of harmonic and subharmonic functions]
{Characterization of harmonic and subharmonic functions via 
mean--value properties}

\author{Pedro Freitas and Jo\~{a}o Palhoto Matos}

\address{Departamento de Matem\'{a}tica, Instituto Superior 
T\'{e}cnico,
Av.Rovisco Pais, 1049-001 Lisboa, Portugal}
\email{pfreitas@math.ist.utl.pt, jmatos@math.ist.utl.pt}
\thanks{The authors were partially supported by FCT, Portugal, 
through program POCTI}



\date{\today}

\dedicatory{}


\begin{abstract}
We give a characterization of harmonic and subharmonic functions in 
terms of their mean values in balls and on spheres. This includes the 
converse of an inequality of Beardon's for subharmonic functions.
We also obtain integral inequalities of Harnack type between these two
means in general domains.
\end{abstract}

\maketitle

\section{Introduction and results}
Let $h$ be a harmonic function in the closed ball centred at a point 
$a$ and 
with radius $r$, $\overline{B}(a,r)$. Then it is well--known that $h$ 
satisfies the mean--value properties
\begin{equation}
    \label{mv1} 
h(a) = \va_{a}(r) := \frac{\ds 1}{\ds |B(a,r)|}\int_{B(a,r)}h(x)\,dx
\end{equation}
and
\begin{equation}
    \label{mv2}
h(a)=\sa_{a}(r) := \frac{\ds 1}{\ds |\partial B(a,r)|}\int_{\partial 
B(a,r)}h(x)ds.
\end{equation}
In both cases, there exist converse results stating that if the 
equality is satisfied, then $h$ must be harmonic. In the first case, 
for
instance, if we assume that $h$ is locally integrable 
on a domain $\Omega$ and satisfies~(\ref{mv1}) whenever 
$\overline{B}(a,r)$
is contained in $\Omega$, then we have that $h$ is harmonic in 
$\Omega$. A similar result holds for~(\ref{mv2}) provided we assume 
now that $h$ is continuous in $\Omega$ --- for these classical results,
see~\cite{abr,gitr}, for instance.

A second, lesser known, set of converse properties to these mean 
value results
assumes that the equalities are satisfied in a set $\Omega$ for all 
integrable harmonic functions defined in $\Omega$. Then, it is 
possible to show that this 
implies that $\Omega$ is a ball~\cite{argo,ben,epsc,kos,kur,paysch}. In
particular,~\cite{ben},~\cite{kos} and~\cite{paysch} consider the 
case where 
the {\it surface} and {\it volume} averages are equal.

Motivated by this, and since we could not find any reference in the 
literature to a converse result of the first type when it is assumed 
that the two averages coincide, in this paper we consider several 
relations between the integral means of harmonic and subharmonic 
functions. We begin by considering the question 
of what can be said when $h$ satisfies
\begin{equation}
    \label{mv3}
    \va_{a}(r) = \sa_{a}(r)
\end{equation}
for balls in a set $\Omega$. We prove that, under suitable 
assumptions, this last equality also implies that $h$ is harmonic. 
More precisely, we have the following
\begin{thm}\label{th} Let $h$ be a continuous function in $\Omega$. 
Then $h$ is harmonic in $\Omega$ if and only if~(\ref{mv3}) 
    is satisfied for all balls whose closure is contained in 
    $\Omega$.
\end{thm}
\begin{rem}{\rm Note that in the equality~(\ref{mv3}) it is not 
assumed a priori that the value of the averages remains fixed when 
$a$ is a 
given point but the radius changes.}
\end{rem}

\begin{rem}{\rm Continuity is essential, since otherwise any function 
coinciding with a harmonic function on all but a set of finite 
points, 
for instance, would clearly satisfy~(\ref{mv3}) but not be harmonic.}
\end{rem}

From the results in~\cite{ben},~\cite{kos} and~\cite{paysch} we have 
that the equality of the averages~(\ref{mv3}) for all harmonic 
functions is possible only for balls. M.
Rao has shown that in fact the averages of harmonic functions on 
general domains must satisfy a 
one--sided inequality of Harnack type~\cite{rao}. Here we improve on 
this result 
establishing also a lower bound by using the methods of~\cite{ben} 
and~\cite{paysch} to obtain the following
\begin{thm}\label{thomega} Let $\Omega$ be a domain satisfying the 
uniform exterior
sphere condition. Then there exist constants $0< c_{1}\leq 1\leq 
c_{2}<\infty$, depending only on $\Omega$, such
that
\[
\frac{\ds c_{1}}{\ds 
|\partial\Omega|}\int_{\partial\Omega}h\,ds\leq
\frac{\ds 1}{\ds |\Omega|}\int_{\Omega}h\,dx \leq \frac{\ds 
c_{2}}{\ds 
|\partial\Omega|}\int_{\partial\Omega}h\,ds
\]
for all non--negative harmonic functions $h$ defined in $\Omega$ and 
which take 
continuous values on the boundary. If either of the constants $c_{1}$ 
or $c_{2}$ can be taken to be equal to one, then $\Omega$ is a ball.
\end{thm}
The proof of the second part of the result also uses a variation of 
a well--known result of Serrin's~\cite{ser}, which we believe to be 
interesting
in its own right -- see Lemma~\ref{varser} for the details.

Another natural question to ask is whether a result similar to 
Theorem~\ref{th} holds for 
subharmonic functions. Here we call a real--valued function $u$ 
defined 
in $\Omega$ subharmonic, if $u$ is continuous and 
for every ball $B$ whose closure is contained in $\Omega$ and every 
harmonic function $h$ in $B$ satisfying $h\geq u$ on $\partial B$ we 
have that $h\geq u$ in $B$~\cite{gitr}. We have that subharmonic 
functions
satisfy
\begin{equation}
    \label{mvs1}
    \va_{a}(r)\leq\sa_{a}(r),
\end{equation}
so it makes sense to ask if this inequality 
characterizes subharmonic functions. In this respect we have the 
following
\begin{thm}\label{tsh} Let $u$ be a continuous function in $\Omega$.
Then $u$ is subharmonic in $\Omega$ if and only if~(\ref{mvs1}) 
holds for all balls whose closure is contained in 
$\Omega$.
\end{thm}

From the proof of Theorem~\ref{thomega} it follows that the 
second inequality in that result also holds for $C^2$ subharmonic 
functions and this can easily be improved to $C^0$ sub-harmonic 
functions. Also for particular domains or in dimension $1$ one easily 
exhibits sequences of positive subharmonic functions for which the 
volume average is bounded and the surface average becomes unbounded. 
This can be improved for relatively general domains provided the 
geometry of the boundary allows a sufficiently rich supply of 
harmonic subfunctions, e.g., if the domain satisfies the 
\emph{uniform exterior sphere condition}.
\begin{thm}\label{harnsub}
    Let $\Omega$ be as above, and $u$ be a non--negative
    subharmonic function in $\Omega$ and continuous in 
$\overline{\Omega}$.
    Then there exists a constant $c_{2}$ as in Theorem~\ref{thomega} 
such that
    \[
    \frac{\ds 1}{\ds |\Omega|}\int_{\Omega}udx \leq \frac{\ds 
c_{2}}{\ds 
    |\partial\Omega|}\int_{\partial\Omega}u\,ds.
    \]
    On the other hand, there exists a sequence of positive 
subharmonic functions $u$ for which the volume average in $\Omega$ 
remains
bounded, while the surface average is unbounded.
\end{thm}

Finally, we consider a property of subharmonic functions proved by 
Beardon~\cite{bear}, which states that if $u$ is subharmonic in a 
set $\Omega$, then
\begin{equation}
    \label{bearcond}
    \sa_{a}(\ka r)\leq\va_{a}(r)
\end{equation}
for all balls whose closure is contained in $\Omega$ and where
\begin{equation}
    \label{beark}
\ka=\left\{
\begin{array}{ll}
    \frac{\disps 1}{\disps 2}, & n=1\eqskip\\
    e^{-1/2}, & n=2\eqskip\\
    \left(\frac{\disps 2}{\disps n}\right)^{1/(n-2)}, & n\geq 3.
\end{array}
\right.
\end{equation}
We prove that this property again characterizes subharmonic 
functions. 
More precisely, we have
\begin{thm}\label{bearprop}
    Let $u$ be a continuous function defined on a domain $\Omega$
    and assume that condition~(\ref{bearcond}) is 
    satisfied for all balls whose closure is contained in $\Omega$ 
and 
    for some $\ka$ satisfying
    \[
    0<\ka\leq \ka_{1}:=-\frac{\disps n}{\disps 4} + \frac{\disps 
1}{\disps 
    2}\sqrt{\frac{\disps n^{2}}{\disps 4}+2n}.
    \]
    Then $u$ is subharmonic in $\Omega$.
\end{thm}
Since the upper bound for $\ka$ in the theorem is larger than (or 
equal to, in the case where $n$ is one) the 
values of $\ka$ given in~(\ref{beark}), this yields a converse of 
Beardon's property.

\section{Harmonic functions}

\subsection{Proof of Theorem~\ref{th}}
Assume that~(\ref{mv3}) holds for a given fixed point 
$a$ in $\Omega$, and all $r$ in $(0,R)$, for some positive number $R$ 
such that 
$\overline{B}(a,R)$ is contained in $\Omega$. Denoting by 
$\omega_{n}$ 
the volume of the unit ball in $\R^{n}$ we have that
\[
\frac{\ds 1}{\ds \omega_{n}r^{n}}\int_{B(a,r)}h(x)dx =
\frac{\ds 1}{\ds n\omega_{n}r^{n-1}}\int_{\partial B(a,r)}h(x)ds
\]
and thus
\[
n\int_{0}^{r}\left[\int_{\partial B(a,t)}h(x)ds\right]dt = r
\int_{\partial B(a,r)}h(x)ds.
\]
Letting
\[
\varphi(r) = \int_{\partial B(a,r)}h(x)ds,
\]
we have that
\[
n\int_{0}^{r} \varphi(t)dt = r\varphi(r).
\]
Since $h$ is continuous, we can differentiate this with respect to 
$r$ 
which gives that $\varphi$ satisfies the differential equation
\[
r\varphi'(r) + (1-n)\varphi(r) = 0.
\]
Hence $\varphi(r) = cr^{n-1}$ for some constant $c$ and
\[
\frac{\ds c}{\ds n\omega_{n}} = \frac{\ds \phi(r)}{\ds 
n\omega_{n}r^{n-1}} = 
\sa_{a}(r),
\]
from which it follows that, for $a$ and $r$ as above, the average 
$\va_{a}(r)$ is independent of $r$. Since $h$ is continuous, the 
value taken by this 
average must be attained by $h$ for some point in the ball $B(a,r)$. 
As this 
holds for arbitrarilly small $r$, we have that this value must 
coincide with $h(a)$. Applying now the converse result to the 
mean--value equality~(\ref{mv1}) completes the proof of 
Theorem~\ref{th}.

\subsection{The integral Harnack inequality}
Let $v$ be the solution of the equation
    \begin{equation} \label{veq}
    \begin{cases}
    \Delta v + 1 = 0, & x\in\Omega\\
    v = 0, & x\in\partial\Omega.
    \end{cases}
    \end{equation}
We have that
\[
\frac{\ds 1}{\ds |\Omega|} \int_{\Omega}h\,dx = -\frac{\ds 1}
{\ds |\Omega|} \int_{\Omega}h\Delta v\,dx =
\frac{\ds 1}{\ds |\Omega|} \int_{\partial\Omega}h\left(-\frac{\ds 
\partial v}{\ds\partial \nu}\right)ds.
\]
Since there exist constants $c_{1}$ and $c_{2}$ such that
\[
0<c_{1}<-\frac{\ds \partial v}{\ds\partial \nu}<c_{2}<\infty,
\]
the inequalities follow.

If $c_{1}$ can be taken to be equal to $1$, we have that
\begin{equation}\label{ineq1}
\int_{\partial\Omega}h\left( \frac{\ds 1}{\ds |\partial\Omega|}+
\frac{\ds 1}{\ds |\Omega|}\frac{\ds \partial v}{\ds\partial 
\nu}\right)ds\leq 0.
\end{equation}
Let $c$ be a positive real number such that
\[
h_{0} = c+\frac{\ds 1}{\ds |\partial\Omega|}+
\frac{\ds 1}{\ds |\Omega|}\frac{\ds \partial v}{\ds\partial 
\nu}
\]
is positive, and choose now $h$ to be equal to $h_{0}$ on the 
boundary. 
Since
\[
\int_{\partial\Omega} \frac{\ds 1}{\ds |\partial\Omega|}+
\frac{\ds 1}{\ds |\Omega|}\frac{\ds \partial v}{\ds\partial 
\nu}ds=0,
\]
replacing $h$ in~(\ref{ineq1}) gives
\[
\int_{\partial\Omega}\left( \frac{\ds 1}{\ds |\partial\Omega|}+
\frac{\ds 1}{\ds |\Omega|}\frac{\ds \partial v}{\ds\partial 
\nu}\right)^{2}ds\leq 0,
\]
from which it follows that $\partial v/\partial \nu$ is constant on 
the boundary. By Serrin's result~\cite{ser}, $\Omega$ must be a ball.

If $c_{2}$ equals $1$, we now have
\[
\frac{\ds 1}{\ds |\Omega|}\int_{\partial\Omega}h\left(-\frac{\ds 
\partial v}
{\ds\partial \nu}\right)ds\leq\frac{\ds 1}{\ds 
|\partial\Omega|}\int_{\partial\Omega}h\,ds
\]
and taking $h=-\partial v/\partial\nu$ on the boundary yields
\[
\frac{\ds 1}{\ds |\Omega|}\int_{\partial\Omega}\left(\frac{\ds 
\partial v}
{\ds\partial \nu}\right)^{2}ds\leq \frac{\ds 1}{\ds 
|\partial\Omega|}\int_{\partial\Omega}\left(-\frac{\ds \partial v}
{\ds\partial \nu}\right)ds = \frac{\ds |\Omega|}{\ds|\partial\Omega|}.
\]
The result now follows from the following lemma, which is a variation 
on
Serrin's result~\cite{ser}.
\begin{lem}\label{varser} Let $v$ be the solution of 
equation~(\ref{veq}).
    Then
    \[
    \frac{\ds 1}{\ds \vdo}\intdo \left(\frac{\ds 
    \partial v}{\ds \partial \nu}\right)^{2}dS \geq 
    \left(\frac{\ds\vo}{\ds|\partial\Omega|}\right)^{2},
    \]
    with equality if and only if $\Omega$ is a ball.
\end{lem}
\begin{proof} 
    By the Cauchy--Schwarz inequality we have that
    \begin{equation}\label{csi}
    \intdo \left(-\frac{\ds \partial v}{\ds \partial \nu}\right)dS
    \leq \vdo^{1/2}\left[\intdo \left(\frac{\ds 
    \partial v}{\ds \partial \nu}\right)^{2}dS\right]^{1/2}.
    \end{equation}
    Thus
    \[
    \frac{\ds 1}{\ds \vdo}\intdo \left(\frac{\ds \partial v}{\ds 
\partial 
    \nu}\right)^{2}dS\geq \left(
    \frac{\ds 1}{\ds \vdo}\intdo \frac{\ds \partial v}
    {\ds \partial \nu}dS\right)^{2} = \left(
    \frac{\ds 1}{\ds \vdo}\intom \Delta vdx \right)^{2} =
    \left(\frac{\ds \vo}{\ds \vdo}\right)^{2}.
    \]
    On the other hand, equality holds in~(\ref{csi}) if and only if 
$\partial
    v/\partial \nu$ is constant on the boundary and thus, by Serrin's 
    result,
    $\Omega$ must be a ball.
\end{proof}

\section{Subharmonic functions}
\subsection{Proof of Theorems~\ref{tsh} and~\ref{harnsub}}

We begin by proving the following
\begin{lem}
    If $u$ is a continuous function in $\Omega$ such that 
$\va_{a}(r)\leq\sa_{a}(r)$ holds for all balls whose closure is 
contained in 
$\Omega$, then $u$ satisfies the maximum principle in $\Omega$.
\end{lem}
\begin{proof}
Let $M$ be the maximum of $u$ in $\overline{\Omega}$, and assume that this 
is attained at an interior point. Define
\[
\mathcal{A} = \left\{ x\in\Omega:u(x) = M\right\}
\]
From the continuity of $u$ we have that $\mathcal{A}$ is closed in 
$\Omega$ and that $\sa$ is a continuous function of the radius. Let 
$x_{0}$ be a point in $\mathcal{A}$, and take $r_{1}$ such that 
$\overline{B(x_{0},r_{1})}$ is in $\Omega$. We have that either there 
exists
$\overline{r}$ in $(0,r_{1})$ such that $\sa_{x_{0}}(\overline{r})<M$ 
or 
not. In the second case, $\sa_{x_{0}}(r)=M$ for all $r$ in 
$(0,r_{1})$, 
which immediately yields that $x_{0}$ is an interior point of 
$\mathcal{A}$.

In the first case, let
\[
r_{*} = 
\inf\left\{r\in(0,\overline{r}):\sa_{x_{0}}(r)=S_{x_{0}}(\overline{r})
\right\}.
\]
Again by continuity, we have that $r_{*}$ is positive and also that 
$\sa_{x_{0}}(r_{*})=\sa_{x_{0}}(\overline{r})<M$. Hence
\[
\begin{array}{lll}
\omega_{n}r_{*}^{n}\va_{x_{0}}(r_{*}) & = 
&{\ds \int_{B(x_{0},r_{*})}}vdx\formsk\\
& = & n\omega_{n}{\ds \int_{0}^{r_{*}}} 
r^{n-1}\sa_{x_{0}}(r)dr\formsk\\
& > & \omega_{n}r_{*}^{n}\sa_{x_{0}}(r_{*})
\end{array}
\]
and thus $\va_{x_{0}}(r_{*})>\sa_{x_{0}}(r_{*})$, contradicting the 
hypothesis.
\end{proof}

{\it Proof of Theorem~\ref{tsh}:} Assume that $w$ satisfies the 
inequality
$\va_{a}(r)\leq\sa_{a}(r)$ but that it is not subharmonic. Then, 
there exists a ball $B$ and a function $h$ harmonic in $\overline{B}$ 
such that
$\overline{B}$ is contained in 
$\Omega$ and $h\geq w$ on $\partial B$, but $h(x_{0})<w(x_{0})$ for 
some 
$x_{0}$ in $B$. Since the function $w-h$ also 
satisfies the inequality but not the maximum principle,
we have a contradiction and hence $w$ must be subharmonic. \qed

{\it Proof of Theorem~\ref{harnsub}:} Let $v$ be subharmonic. Let $h$ 
be the harmonic solution of the Dirichlet problem with boundary 
values specified by $v$. Then $v-h$ is subharmonic and the maximum 
principle shows that $v\leq h$. Hence Theorem \ref{thomega} applied 
to $h$ implies
\[
\frac{1}{|\Omega|}\int_\Omega v \, dx\leq 
\frac{1}{|\Omega|}\int_\Omega h \, dx\leq 
\frac{c_2}{|\partial\Omega|}\int_{\partial\Omega}h\, ds = 
\frac{c_2}{|\partial\Omega|}\int_{\partial\Omega}v\, ds.
\]

For the second part of the statement consider first the trivial one 
dimensional case. Assume (without loss of generality) $\Omega=(0,1)$ 
and let, for $k$ a positive integer,  $v_k(x)=\max\{1, 
k-k^2x,k+k^2(x-1) \}$. Then the surface average is $k$ and the volume 
average is bounded.

To generalize the one dimensional case we will need for each element 
of the sequence a maximum of a finite number of harmonic functions 
which is large in a neighborhood of $\partial\Omega$ and small 
elsewhere. To deal with a fairly general geometry of $\partial\Omega$ 
requires using in the construction harmonic functions other than 
affine functions and will require a covering argument to control the 
oscillation of these on $\partial \Omega$. For domains satisfying the 
uniform exterior sphere condition translating and rescaling 
fundamental solutions will be enough. Details follow below.

Let $k>2$ be an integer and $n>1$. For each $y\in\partial\Omega$, let 
$\overline{B}_{\delta_k}(\overline{y}^k)$ denote a closed ball such 
that 
$\overline{B}_{\delta_k}(\overline{y}^k)\cap\overline{\Omega}=\{y\}$ 
and $\delta_k>0$. Here $(\delta_k)$ is a sequence of positive numbers 
smaller than the radius in the uniform exterior sphere condition and 
verifying $\delta_k\to 0$.

Denote the fundamental solution of the Laplace operator by 
$\Psi(|x|)$.
Define, for each $k>2$ and $y\in\partial \Omega$
\[
v_y^k(x)=-k^2\Psi(|x-\overline{y}^k|/\delta_k)+k^2\Psi(|y-\overline{y}^k|/\delta_k)+k.
\]

Clearly $v_y^k(y)=k$ and $v_y^k<k$ on $\partial\Omega\setminus\{y\}$. 
By continuity and compactness there is a finite number (dependent on 
$k$) of points $y_j^k\in\partial\Omega$ and corresponding balls 
$B_{r_{jk}}(y_j^k)$ forming a covering of $\partial\Omega$ such that 
each $v_{y_j}^k>k-1$ in $B_{r_{jk}}(y_j^k)$.

Let $v^k=\max\{\max_j\{v_{y_j}^k\},1\}$. Then $v^k$ is subharmonic, 
$k\geq v^k\geq 1$ in $\overline{\Omega}$ and $v^k>k-1$ on $\partial 
\Omega$.

Let $\Omega_\delta=\{x\in\Omega:\dist(x,\partial\Omega)\geq 
\delta\}$. If $z\in\Omega_{\delta_k}$ then $|z-\overline{y}_j^k|\geq 
2\delta_k$ for all $\overline{y}_j^k$ and $v^k_{y_j}(z)\leq 
-k^2\Psi(2)+k^2\Psi(1)+k$ implying that, for sufficiently large $k$, 
$v^k=1$ in $\Omega_{\delta_k}$. This allows, for sufficiently large 
$k$,  the estimates
\begin{align*}
\int_{\partial\Omega}v^k\,ds & >|\partial\Omega|(k-1),\\
\int_\Omega v^k\, dx & = \int_{\Omega_{\delta_k}} 1 + 
\int_{\Omega\setminus\Omega_{\delta_k}} v^k\, dx \leq |\Omega| + k 
|\Omega\setminus \Omega_{\delta_k}|.
\end{align*}
The result follows by choosing the sequence $(\delta_k)_{k>2}$ in such a 
way that it implies bundedness of the last expression.
\qed

\subsection{Beardon's property}

The converse of Beardon's property follows from three facts: a 
standard mollification argument that shows the result holds if it 
holds for smooth functions, a computation for smooth functions that 
shows that Beardon's inequalities holding for balls centered at a 
point $x_0$ imply the Laplacian is non--negative at $x_0$, provided 
$\ka<\ka_1$, and a well known device to extend the latter result to 
$\ka\leq \ka_1$. We present these as three lemmas.


\begin{lem}
Assume $u$ is a continuous function for which (\ref{bearprop}) holds 
and is not subharmonic. Then a smooth subharmonic function exists for 
which the same properties hold.
\end{lem}

\begin{proof}Let $u$ be such a continuous function. Then, in some 
point $\overline{x}$ in the interior of some ball $B\Subset\Omega$, 
$u$ is bigger then the solution $h$ of Dirichlet's problem for the 
the Laplace equation in $B$ having $u$ as boundary data. Denote 
$\gamma\equiv u(\overline{x})-h(\overline{x})>0$.

Let $\delta=\dist(B,\partial \Omega)$ and $0<\epsilon<\delta/4$ and 
consider the usual mollifiers $\rho_\epsilon$ supported in 
$B_\epsilon(0)$, $\varphi$ a continuous cut-off function which is $0$ 
in $\Omega\setminus\Omega_{\delta/4}$ and $1$ in $\Omega_{\delta/2}$ 
and extend $u$ by $0$ in the complement of $\Omega$.

Then Beardon's inequality holds for $\rho_\epsilon\ast (u\varphi)$ in 
$\Omega_{3\delta/4}$. To check this statement let 
$B_r(a)\Subset\Omega_{3\delta/4}$ and notice that
\[
\begin{split}
& \frac{1}{r^{n-1}\ka^{n-1}\omega_n}\int_{\partial 
B_{kr}(a)}\rho_\epsilon\ast(u\varphi)\,dS\\ & \qquad 
=\frac{1}{r^{n-1}\ka^{n-1}\omega_n}\int_{\partial 
B_{kr}(a)}\left(\int_{B_\epsilon(0)}\rho_\epsilon(y)u(x-y)\varphi(x-y)\,dy\right)\,dS(x)\\
& \qquad = 
\frac{1}{r^{n-1}\ka^{n-1}\omega_n}\int_{B_\epsilon(0)}\left(\rho_\epsilon(y)\int_{\partial 
B_{kr}(a)}u(x-y)\,dS(x)\right)\,dy\\
& \qquad \leq 
\frac{n}{r^{n}\omega_n}\int_{B_\epsilon(0)}\left(\rho_\epsilon(y)\int_{B_{r}(a)}u(x-y)\,dx\right)\,dy 
\\
& \qquad = 
\frac{n}{r^{n}\omega_n}\int_{B_{r}(a)}\left(\int_{B_\epsilon(0)}\rho_\epsilon(y)u(x-y)\,dy\right)\,dx\\
& \qquad = 
\frac{n}{r^{n}\omega_n}\int_{B_{r}(a)}\rho_\epsilon\ast(u\varphi)\, 
dx.
\end{split}
\]

As $\rho_\epsilon\ast (u\varphi)\to u$ uniformly on $\partial B\cup 
\{\overline{x}\}$ as $\epsilon\to 0$, we have, for sufficiently small 
$\epsilon>0$, $\rho_\epsilon\ast(u\varphi)-\gamma/2<h$ on $\partial 
B$ and  
$\rho_\epsilon\ast(u\varphi)(\overline{x})-\gamma/2>h(\overline{x})$.

Hence $\rho_\epsilon\ast (u\varphi)-\gamma/2$ has all the desired 
properties in $\Omega_{3\delta/4}$.
\end{proof}

\begin{lem}\label{lem9}
Let $u\in C^2(B_R(a))$ and assume (\ref{bearprop}) holds for $u$ for 
some $\ka$ satisfying $0<\ka<\ka_1$ and  all $r$ satisfying $0<r<R$. 
Then $\Delta u(a)\geq 0$.
\end{lem}

\begin{proof}Assume $a=0$. Beardon's inequality can be rewritten as
\[
\begin{split}
0 & \leq \frac{\disps n}{\disps r^{n}}{\disps \int_{0}^{r}}\left(
{\disps \int_{\partial 
B_{t}(0)}}u(x)dS(x)\right)dt - \frac{\disps 1}{\disps 
r^{n-1}\ka^{n-1}}
{\disps \int_{\partial B_{\ka r}(0)}}u(x)\,dS(x)\\
 & =  \frac{\disps n}{\disps r^{n}}{\disps \int_{0}^{r}}\left(t^{n-1}
 {\disps \int_{\partial B_{1}(0)}}u(tx)dS(x)\right)\,dt-
 {\disps\int_{\partial B_{1}(0)}}u(\ka rx)\,dS(x).
\end{split}
\]

Let now
\[
\psi(r) = n{\disps \int_{0}^{r}}\left(t^{n-1}
 {\disps \int_{\partial B_{1}(0)}}u(tx)dS(x)\right)\,dt-
 r^{n}{\disps\int_{\partial B_{1}(0)}}u(\ka rx)\,dS(x).
\]
As $\psi\geq 0$ we have that
\[
\varliminf_{r\to 0}\frac{\disps \psi(r)}{\disps r^{n+2}} \geq 0,
\]
which in turn equals
\[
\lim_{r\to 0}\frac{\disps \psi'(r)}{\disps (n+2)r^{n+1}},
\]
provided the latter limit exists. Now
\begin{equation*}
\begin{split}
\psi'(r)  & =  nr^{n-1}{\disps \int_{\partial 
B_{1}(0)}}u(tx)\,dS(x)-
 nr^{n-1}{\disps \int_{\partial B_{1}(0)}}u(\ka rx)\, dS(x)\\
 & \quad - r^{n}{\disps \int_{\partial B_{1}(0)}}\ka \nabla u(\ka r 
 x)\cdot x\,dS(x)\\
 & = nr^{n-1}{\disps \int_{\partial B_{1}(0)}}\left[{\disps \int_{\ka 
 r}^{r}} \frac{\disps d}{\disps dt}\left( u(tx)\right)dt\right]
 dS(x)-r^{n}{\disps \int_{\partial B_{1}(0)}}\ka \nabla u(\ka r 
 x)\cdot x\,dS(x)\\
  & =  nr^{n-1}{\disps \int_{\partial B_{1}(0)}}\left[{\disps 
\int_{\ka 
 r}^{r}}\nabla u(t x)\cdot x dt\right]dS(x) -r^{n}{\disps 
\int_{\partial B_{1}(0)}}
 \ka \nabla u(\ka r x)\cdot x\,dS(x)\\
&  = nr^{n-1}{\disps \int_{\ka r}^{r}}\left[
 {\disps \int_{\partial B_{1}(0)}}\nabla u(t x)\cdot x\,dS(x)\right]dt
 -r^{n}{\disps \int_{\partial B_{1}(0)}}\ka \nabla u(\ka r x)\cdot 
x\,dS(x)\eqskip\\
&  =   nr^{n-1}{\disps \int_{\ka r}^{r}}\left[
 {\disps \int_{B_{1}(0)}}t\Delta u(t x) dx\right]dt-
 r^{n+1}{\disps \int_{B_{1}(0)}}\ka^{2} \Delta u(\ka r x)dx,
\end{split}
\end{equation*}
where the last step follows by applying the divergence theorem. Thus,
\[
\begin{split}
{\disps \lim_{r\to 0}}\frac{\disps \psi'(r)}{\disps r^{n+1}} & = 
{\disps \lim_{r\to 0}}\left[\frac{\disps n}{\disps r^{2}}{\disps 
\int_{\ka r}^{r}}\left[
 {\disps \int_{B_{1}(0)}}t\Delta u(t x) dx\right]dt-
 \ka^{2} {\disps \int_{B_{1}(0)}}\Delta u(\ka r x)dx\right]\eqskip\\
 & =  {\disps \lim_{r\to 0}}\left[\frac{\disps n}{\disps 
r^{2}}{\disps \int_{\ka r}^{r}}\left[
 {\disps \int_{B_{1}(0)}}t\Delta u(t x) dx\right]dt\right]-
 {\disps \ka^{2} \int_{B_{1}(0)}}\Delta u(0)dx.
\end{split}
\]
Since
\[
{\disps \lim_{r\to 0}} \left[\frac{\disps n}{\disps r^{2}}{\disps 
\int_{\ka r}^{r}}\left[
 {\disps \int_{B_{1}(0)}}t\Delta u(t x) dx\right]dt\right] = 
 \frac{\disps n}{\disps 2}
{\disps \int_{B_{1}(0)}}\left(1 -\ka\right) \Delta u(0) dx,
\]
we finally obtain that
\[
\left[ \frac{\disps n}{\disps 2}(1-\ka)-\ka^{2}\right]\Delta 
u(0)\geq0,
\]
and thus if $\ka\in(0,\ka_{1})$, it follows that $\Delta u(0)\geq 0$.
\end{proof}

The next lemma is only relevant in establishing the converse of 
Beardon's property in the case $n=1$.
\begin{lem}
The previous lemma holds if $\ka=\ka_1$.
\end{lem}
\begin{proof}
Again assume $a=0$. Assume Beardon's inequality holds for some smooth 
function $u$ with $\ka=\ka_1$. Let $v(x)=|x|^2$. Then
\[
\frac{1}{\omega_n|\ka|^{n-1}r^{n-1}}\int_{\partial B_{\ka 
r}}v(x)\,dS=\ka^{2} r^{2}
\]
and
\[
\frac{n}{\omega_n r^{n}}\int_{B_{r}}v(x)\,dx 
=\int_0^r\left(\int_{\partial B_{\rho}}|\rho|^2\,dS\right)d\rho = 
\frac{n}{n+2}r^{2}.
\]
Hence Beardon's property holds for $v$ with $\ka=\ka_0$, hence it 
holds for $u+\epsilon v$ for any $\epsilon>0$ and with $\ka=\ka_0$ 
and consequently $\Delta u(0)+2n\epsilon = \Delta(u+\epsilon 
v)(0)\geq 0$ for any $\epsilon>0$. But then $\Delta u(0)\geq 0$.
\end{proof}

From the proof of Lemma~\ref{lem9} we have that if 
$\va_{a}(r)\leq\sa_{a}(\ka r)$ for $\ka>\ka_{1}$ then
\[
\left[ \frac{\disps n}{\disps 2}(1-\ka)-\ka^{2}\right]\Delta 
u(a)\leq0,
\]
and we again conclude that $\Delta u(a)\geq0$. This of course already 
followed from Theorem~\ref{tsh} (of which this argument gives a 
different proof). On the other hand, we have that for $u(x) = |x|^{2p}$
($p$ integer)
\[
\va_{0}(r) = \frac{\disps n r^{2p}}{\disps 2p+n},
\mbox{ while }
\sa_{0}(r) = \ka^{2p}r^{2p},
\]
giving that $\va_{0}(r)\leq\sa_{0}(\ka r)$ provided
\[
\ka\geq \left(\frac{\disps n}{2p+n}\right)^{1/(2p)},
\]
which converges to one as $p\to\infty$. Thus we see that the best possible value
of $\ka$ for this property of subharmonic functions is one.




\end{document}